\documentclass [12pt]{article}
\usepackage{amsmath}
\usepackage{amsfonts, amssymb}

\newtheorem{thm}{Theorem}
\newtheorem{lemma}[thm]{Lemma}

\def\R{{\mathbf R}}
\def\E{{\mathbf E}}
\def\P{{\mathbf P}}
\def\Z{{\mathbf Z}}
\def\N{{\mathbf N}}

\def\be#1\ee{\begin{equation}#1\end{equation}}
\newcommand{\bea}{\begin{eqnarray}}
\newcommand{\eea}{\end{eqnarray}}
\newcommand{\beaa}{\begin{eqnarray*}}
\newcommand{\eeaa}{\end{eqnarray*}}

\def\EE{{\mathcal E}}
\def\F {{\mathcal{F}}}
\def\ga{{\gamma}}

\begin{document}

\title{Cyclic behavior of maxima for sums of independent variables
\footnote{Work supported by grants RFBR 13-01-00172, 11-01-12104-ofi\_m,
FFP 2010-1.1.-111-128-033, and by SPSU research grant 6.38.672.2013.}
}

\author
{M. A. Lifshits}
\date{ }
\maketitle

\begin{abstract}
In a recent author's work the cyclic behavior of maxima in a hierarchical 
summation scheme was discovered. In the present note we show how the same
phenomenon appears in the scheme of conventional summation: the distribution
of maximum of $2^n$ independent copies of a sum of $n$ i.i.d. random
variables  approaches, as $n$ grows, some helix in the space of distributions. 
\end{abstract}

\section{Introduction and main results}
\label{s:notree}

In a recent author's work \cite{Lif} he proved a cyclic limit theorem for a 
hierarchical summation scheme that also admits an interpretation in terms of
a branching random walk. It turned out that for a branching random walk with two 
descendants and symmetric Bernoulli displacement of particles the distribution
of location of the extremal (maximal or minimal) particle moves, as the number of
generation grows, along a helix in the space of distributions. If one takes shifts into 
account, then the distributions of positions of extremal particles rotate along
a closed curve.     
Moreover, since the motion of the sequence  goes along the limiting helix more and
more slowly, all points of the helix are limiting points under appropriate shift 
normalization; there is no unique limiting distribution for the positions of extremal
particles. 

The cyclic effect, implicitly and without geometric interpretation, is also contained
in Theorem 1 of classical Bramson work \cite{Br1}. The advantage  
of his results is in consideration of more general branching mechanism and handling a.s.
convergence instead of convergence in distribution. 

In the present note we show how the same cyclic effect may be obtained by conventional 
summation: the distribution of maximum of $2^n$ independent copies 
of a sum of $n$ independent identically distributed random variables are attracted,
as $n$ grows, to a helix in the space of distributions. 

Let $(\xi_i)_{i\in \N}$ be {\it integer-valued}\ i.i.d. random variables.   
Consider the sum $S_n:=\sum_{i=1}^n \xi_i$, and let 
$S_n^{(j)}$, $1\le j\le 2^n$,  be independent copies of $S_n$. 
We are interested in the behavior of variables $M_n:=\max_{j\le 2^n} S_n^{(j)}$. 

In what concerns the distribution of initial random variables, we assume that 
\be \label{moments}
    \E|\xi_1|<\infty\quad \textrm{and}\quad \E\exp\{\ga \xi_1\}<\infty,\ \forall \ga>0.
\ee

Let $\omega$ denote the upper bound of the distribution, i.e. 
\[ 
   \omega:= \sup \{m\in \N: \P(\xi_1=m)>0 \}. 
\]
Assume that one of the following two conditions holds: either

$(i)$ \qquad $\omega=\infty$,

\noindent or

$(ii)$ \qquad  $\omega<\infty$\ and\ $\P\left(\xi_1=\omega\right)<1/2$. 
\medskip

The sense of these conditions is explained in the following lemma on the behavior
of the cumulant $L(\ga):=\ln \E \exp\{\ga \xi_1\}$.

\begin{lemma} \label{l1}
Let \eqref{moments} be satisfied. Then condition $\omega=\infty$ implies
\[ 
    \lim_{\ga\to +\infty} [ L(\ga)-\ga L'(\ga)] = -\infty,  
\]
while condition $\omega<\infty$ implies
\[ 
    \lim_{\ga\to +\infty} [ L(\ga)-\ga L'(\ga)] = \ln\P(\xi_1=\omega).
\]
\end{lemma}

Since the cumulant $L(\cdot)$ is convex, the function $L(\ga)-\ga L'(\ga)$ is
non-decreasing. Moreover, it is continuous and vanishes at $\gamma=0$.  
If we assume in addition that condition 
\eqref{moments} and one of the conditions $(i)$ or $(ii)$ hold, then 
it follows from Lemma \ref{l1}
that
\[ 
    \lim_{\ga\to +\infty} [ L(\ga)-\ga L'(\ga)] < \ln (1/2).
\]
Therefore, a solution of equation
\be \label{gammastar}
    L(\ga)-\ga L'(\ga) = \ln (1/2)
\ee 
on $(0,+\infty)$ exists. Let denote it $\ga_*$ and let
$\rho_*:= L'(\ga_*)$. Notice also that, whenever $(i)$ or
$(ii)$ holds, the common distribution of r.v. $\xi_i$ is non-degenerated 
(i.e. it is not concentrated in a single point),
hence the solution of equation  \eqref{gammastar} is unique.

\begin{thm} \label{cycle_notree} Assume that condition 
\eqref{moments} and one of conditions $(i)$  or $(ii)$ hold. Let  
$\rho_*,\gamma_*$ be defined by equation $(\ref{gammastar})$. Then 
\be \label{bound_gen}
   \P\left\{ M_n < \rho_* n - \frac{\ln n}{2 \gamma_*} +z \right\}
   =\exp\left\{ - \, \frac{\exp\{-\gamma_* z\} (1+o(1))}
    {\sqrt{2\pi}\sigma(\ga_*)(1-e^{-\ga_*})}  \right\},
\ee
where $\sigma(\cdot)^2=L''(\cdot)$, uniformly over
\footnote{In other words, we consider $z$ such that the expression
in the left hand side is an integer.}  
\[
    z\in I\bigcap  \left[\Z - \rho_* n  + \frac{\ln n}{2 \gamma_*} \right]
\]   
for any bounded interval $I$.  
\end{thm}

We may rewrite \eqref{bound_gen} as
\be \label{bound_geni}
   \P\left\{ M_n < m \right\}
   =\exp\left\{ - \exp\{-\gamma_* (m-a_n)\} (1+o(1)) \right\},
   \qquad m\in \Z,
\ee
where
\[
  a_n:= \rho_* n   
        -\frac{\ln[\sqrt{2\pi n}\sigma(\ga_*)(1-e^{-\ga_*})]}{\ga_*}\ .
\]

For any $a\in\R$ let $\F^a$ denote the distribution on integer numbers given by
the relation
\[
   \F^a((m,+\infty))
   = \exp\left\{ - \exp\{-\gamma_* (m-a)\}\right\}, \qquad m\in\Z.
\]
Then $(\F^a)_{a\in\R}$ is a curve in the space of distributions. It is
natural to interpret it as a helix because there is a 1-periodicity
up to a shift: $\F^{a+1}\{m+1\}=\F^{a}\{m\}$. 
Equation \eqref{bound_geni} shows that the distribution of r.v. $M_n$ 
is uniformly approximated by an element $\F^{a_n}$ of this helix, and
after appropriate centering it is approximated by an element
$\F^{[a_n]}$ of the helix turn $(\F^a)_{0\le a <1}$.
Moreover, any distribution $(\F^a)_{0\le a < 1}$ is a limit of some subsequence
of centered distributions of $M_n$.

The proof of Theorem \ref{cycle_notree} given in the next section is based
upon a theorem on large deviations due to V.V. Petrov.

Clearly, the number of sums $2^n$ may be replaced, with obvious changes,
by any other exponentially growing sequence.
\bigskip

Let us consider Bernoulli case as an example. Let
$\xi_i=B_i$ be independent random variables following non-symmetric
Bernoulli distribution, i.e.
\[
   \P(B_i=1)=1-\P(B_i=-1)= p < 1/2 \ .
 \]
Let $q=1-p$ and define the shift coefficient $\rho_*$ by equation 
\be \label{Bdrift}
     2p^\rho q^{1-\rho}= \rho^\rho(1-\rho)^{1-\rho}.
\ee 
  
We also need two additional constants:
$\kappa:=\frac{p(1-\rho_*)}{q\rho_*}\in (0,1)$ and $\beta:=2\pi\rho_*(1-\rho_*)$.
Then the result of Theorem \ref{cycle_notree} takes the following form.
 
 \begin{thm} \label{cycle_notree_b} 
 The representation
 \be \label{bound_ber}
     \P\left\{ M_n< \rho_* n -\frac{\ln(\beta n)}{2|\ln\kappa|} +z \right\}
     =\exp\left\{ - \frac{\kappa^{z}}{1-\kappa}\, (1+o(1))\right\},
 \ee
 holds uniformly over
 \[
    z\in I\bigcap  \left[\Z - \rho_* n +\frac{\ln(\beta n)}{2|\ln\kappa|}\right]
 \]   
 for any bounded interval $I$.
 \end{thm}

{\bf Remark.} If $p \ge \tfrac 12$, then neither of conditions $(i), (ii)$  
holds. Equation \eqref{gammastar} has no solutions, thus 
Theorem \ref{cycle_notree} does not apply. 

\section{Proofs}
\label{s:proofs}

{\bf Proof of Lemma \ref{l1}.} \ 

1) In view of convexity of function $L(\cdot)$ and by $L(0)=0$, the function 
$\ga\to \ga L'(\ga)- L(\ga)$ is increasing and non-negative. Therefore, 
there exists a non-negative limit  
\[ 
   C:= \lim_{\ga\to +\infty} [ \ga L'(\ga)- L(\ga)].
\]
Assume that $C<\infty$. Then for any $\ga>0$ we have
$\ga L'(\ga)- L(\ga)\le C$.
By integrating inequality
$\frac{L'(\ga)}{L(\ga)+C}\le \frac 1\ga$\, ,  we obtain for any $\ga>1$
\[
 \ln(L(\ga)+C)-\ln(L(1)+C) \le \ln \ga,
\]
hence,
\be \label{sublinear}
  L(\ga)\le (L(1)+C)\ga- C.
\ee
On the other hand, for any $m\in \N$ we have
\[
   L(\ga) \ge \ln \P(\xi_1=m) +\ga m.
\]
Therefore, condition $\omega=\infty$ and \eqref{sublinear} are incompatible. 
It follows that $\omega=\infty$ yields $C=+\infty$.

2) Let $\omega<\infty$. Then, without loss of generality, 
we may assume that $\omega=0$, since the change of variable $\xi$
by r.v. $\xi-\omega$ does not affect the function $L(\ga)-\ga L'(\ga)$. 
Furthermore, we have
\begin{eqnarray*} 
   \lim_{\ga\to +\infty} L(\ga) 
   &=&  \lim_{\ga\to +\infty} \ln\left( \P(\xi_1=0)
       +\sum_{j=1}^\infty \P(\xi_1=-j)e^{-\ga j}\right)
 \\
 &=& \ln\P(\xi_1=0).
\end{eqnarray*}
On the other hand,
\[
 \ga L'(\ga)= -\ga \left(\E e^{\ga \xi_1}\right)^{-1} 
             \sum_{j=1}^\infty \P(\xi_1=-j) j e^{-\ga j},
\]
thus for $\ga\to\infty$ we have
\[
 |\ga L'(\ga)| \le \ga e^{-\ga} \left(\P(\xi_1=0) \right)^{-1} 
             \sum_{j=1}^\infty j e^{-\ga (j-1)} \to 0.
\]
Therefore,
\[
  \lim_{\ga\to +\infty} [ L(\ga)- \ga L'(\ga)] 
  = \lim_{\ga\to +\infty} L(\ga) - \lim_{\ga\to +\infty} \ga L'(\ga) 
  =   \ln\P(\xi_1=0).
 \]
$\Box$

{\bf Proof of Theorem \ref{cycle_notree}.} \ 
It is obvious that
\begin{eqnarray}
\nonumber
\ln  \P\left\{ M_n < \rho_* n - \frac{\ln n}{2 \gamma_*} +z \right\}
&=& 2^n \ln  \P\left\{ S_n < \rho_* n - \frac{\ln n}{2 \gamma_*} +z \right\}
\\ \nonumber
&=& 2^n \ln  \left(1-\P\left\{ S_n \ge \rho_* n - \frac{\ln n}{2 \gamma_*} +z \right\}\right)
\\ \label{e2n}
&=& - 2^n \ \P\left\{ S_n \ge \rho_* n - \frac{\ln n}{2 \gamma_*} +z \right\} (1+o(1)).
\end{eqnarray}
We evaluate the latter large deviation probability by using a V.V. Petrov theorem,
see \cite{Petr_TVP} or Complement 2 in \cite[\S4, Chapter VIII]{Petr}. 
According to this theorem, it is true that
\be \label{lda}
  \P\left\{ S_n \ge nx \right\} 
  = \frac{\exp(-n(\ga x- L(\ga)))}{\sqrt{2\pi L''(\gamma) n}\, (1- e^{-\ga})} \ 
    \left( 1+O\left(\frac 1n\right)\right),
\ee
where $\ga$ is the unique solution of equation $L'(\gamma)=x$.

In our case
\[ 
   x=x_n:= \rho_*  - \frac{\ln n}{2 \gamma_*n} + \frac{z}{n} \to \rho_*. 
\]
Therefore, the denominator of fraction \eqref{lda} is equivalent to 
$\sqrt{2\pi L''(\ga_*) n}\, (1- e^{-\ga_*})$, where $\ga_*$ denotes the solution 
of equation $L'(\gamma)=\rho_*$, which is (according to the definition of $\rho_*$) 
the same $\ga_*$ that emerges in the assertion of Theorem \ref{cycle_notree}.
Therefore, it remains to study the expression
\begin{eqnarray}
    \nonumber
    \EE&:=& \frac{\exp(-n(\ga x- L(\ga)))}{\sqrt{n}}
       := \exp(-n(\ga_n x_n- L(\ga_n))-\tfrac{\ln n}{2})
\\   \label{ee}
    &=& \exp(-n(\ga_n x_n  \pm \ga_* \rho_* - L(\ga_n) \pm L(\ga_*))-\tfrac{\ln n}{2}).
\end{eqnarray}
Taylor expansions yield
\begin{eqnarray*}
    \ga_n x_n -\ga_* \rho_* 
    &=& \rho_*(\ga_n-\ga_*)+\ga_*(x_n-\rho_*)+(\ga_n-\ga_*)(x_n-\rho_*)
\\ 
    &=&  \rho_*(\ga_n-\ga_*) +\ga_*\left( \frac zn-\frac{\ln n}{2\ga_*n}\right) 
         +O\left(\frac{(\ln n)^2}{n^2}\right)
 \\ 
    &=&  \rho_*(\ga_n-\ga_*) + \frac {\ga_*z}n-\frac{\ln n}{2n} 
         +O\left(\frac{(\ln n)^2}{n^2}\right);
\\
    L(\ga_n)-L(\ga_*) &=& L'(\ga_*)(\ga_n-\ga_*)+ O\left((\ga_n-\ga_*)^2\right)
\\
     &=& \rho_*(\ga_n-\ga_*)+ O\left((x_n-\rho_*)^2\right)
\\
     &=& \rho_*(\ga_n-\ga_*) +O\left(\frac{(\ln n)^2}{n^2}\right).
\end{eqnarray*}
By plugging these expressions into \eqref{ee} and by using equation \eqref{gammastar}, 
we obtain
\[
  \EE = \exp(-\ga_* z-n\ga_*\rho_*+  n L(\ga_*))\ (1+o(1))
  = 2^{-n} \exp(-\ga_* z)\ (1+o(1)).
\]
By collecting together the obtained estimates, we get
\[
 \P\left\{ S_n \ge \rho_* n - \frac{\ln n}{2 \gamma_*} +z \right\}
 =  2^{-n}\, \frac{\exp(-\ga_* z)}{ \sqrt{2\pi L''(\ga_*)}\, (1- e^{-\ga_*})}   \ (1+o(1)).
\]
By definition, $\sqrt{L''(\ga_*)}=\sigma(\gamma_*)$. By plugging the 
obtained expression into \eqref{e2n}, we obtain the desired expression 
\eqref{bound_gen}. The uniformity of the bound also follows from the 
above mentioned V.V.Petrov theorem.\ $\Box$
\bigskip 

{\bf Proof of Theorem \ref{cycle_notree_b}.} \ For Bernoulli distribution
we have $L(\ga)= \ln(pe^\ga+q)$ and $\rho:=L'(\ga)=\tfrac{pe^\ga}{pe^\ga+q}$.
By representing $\ga$ via $\rho$, we obtain 
\be \label{ega}
   e^\ga=\tfrac{\rho q}{p(1-\rho)}. 
\ee
Equation \eqref{gammastar} for Bernoulli variables takes the form
\[ 
   \ln(pe^\ga+q) - \ga \rho = -\ln 2.
\]     
By plugging in the expressions for $\ga$ and $e^\ga$, we obtain
\[
  \ln\left(\frac{q}{1-\rho}\right) 
  - \rho \ln\left(\frac{\rho q}{p(1-\rho)}\right) =  -\ln 2 
\]
and arrive at equation \eqref{Bdrift}.
Furthermore, by using definitions  of $\sigma(\cdot)$ and \eqref{ega}, 
we find
\begin{eqnarray*}
\sigma^2(\ga)&=& L''(\ga)= \frac{d}{d\ga}\ \frac{p e^\ga}{pe^\ga+q}
\\
&=&  \frac{p q e^\ga}{(pe^\ga+q)^2} = \rho(1-\rho),
\end{eqnarray*}
whereas $\sigma(\ga_*)= [\rho_*(1-\rho_*)]^{1/2}$.

For an arbitrary $z$, let $z_1:= z-\frac{\ln\beta}{2|\ln\kappa|}$ and consider
the corresponding bound for the maximum in the left hand side of
\eqref{bound_gen}. Since \eqref{ega} implies
\[ 
   \ga_*= \ln\left(  \frac{\rho_* q}{p(1-\rho_*)} \right)  = |\ln\kappa|, 
\]
we have
\begin{eqnarray*}
\rho_* n - \frac{\ln n}{2 \gamma_*} +z_1  
&=& \rho_* n - \frac{\ln n}{2 |\ln \kappa|} +z -\frac{\ln\beta}{2|\ln\kappa|}
\\
 &=& \rho_* n - \frac{\ln (\beta n)}{2 |\ln \kappa|} +z,
\end{eqnarray*}
thus we have arrived to the similar bound in the left hand side of 
\eqref{bound_ber}.

It remains to compare the right hand sides of \eqref{bound_gen} and 
\eqref{bound_ber}. 
Indeed, we have
\begin{eqnarray*}
      \frac{\exp\{-\gamma_* z_1\}} {\sqrt{2\pi}\sigma(\ga_*)(1-e^{-\ga_*})} 
      &=&
      \frac{\kappa^{z_1}} {\sqrt{2\pi\rho_*(1-\rho_*)}(1-\kappa)}
\\
       &=&
      \frac{\kappa^{z-\frac{\ln\beta}{2|\ln\kappa|}}} {\sqrt{\beta}(1-\kappa)} 
      =
       \frac{\kappa^{z}\sqrt{\beta}}{\sqrt{\beta}(1-\kappa)} 
      =   \frac{\kappa^{z}}{1-\kappa}\ ,
\end{eqnarray*}    
thus the right hand side parts of \eqref{bound_gen} and \eqref{bound_ber} 
coincide  and Theorem \ref{cycle_notree_b} follows from Theorem \ref{cycle_notree}.
$\Box$
\medskip

\bibliographystyle{amsplain}

\bigskip

\bigskip

St.Petersburg State University and MAI, Link\"oping University.

email: {\tt  lifts@mail.rcom.ru}

\end{document}